\documentclass[10pt]{article}

\usepackage{a4wide}
\usepackage{amssymb}
\usepackage{amsfonts}
\usepackage{amsmath}
\input xy
\xyoption{arrow} \xyoption{matrix}

\date{}

\newtheorem{proposition}{Proposition}[section]
\newtheorem{theorem}[proposition]{Theorem}
\newtheorem{lemma}[proposition]{Lemma}

\newtheorem{corollary}[proposition]{Corollary}

\def\der{\partial }

\def\nFM0{{\nu }_{F,M_0}}
\def\nFN0{{\nu }_{F,N_0}}
\def\nGN0{{\nu }_{G,N_0}}

\def\N0{ {\bf N}_0 }

\def\ra{\rightarrow}

\def\Xpm{X^{\pm }}

\def\s{\sigma}

\def\l1{{\lambda}_1}

\def\a{\alpha}
\def\a0{ {\alpha }_0}
\def\a1{ {\alpha }_1}

\def\l{\lambda}

\def\nFGM0{{\nu }_{F,G,M_0}}

\def\nFN0{{\nu}_{F,N_0}}

\def\sm{{\sigma}^m}

\def\sm1{{\sigma}^{-1}}

\def\smtp1{{\sigma}^{-t+1}}

\def\S1{S^{-1}}

\def\Xpm1{X^{\pm 1}_1}

\def\sPM1{{\sigma }^{\pm 1}}
\def\sMP1{{\sigma }^{\mp 1 }}

\def\di{{\rm d.ind}}

\def\L{\Lambda}

\def\Ytm1{Y^{t-1}}
\def\Yim1{Y^{i-1}}

\def\CL{{\cal L}}

\def\CN{{\cal N}}

\def\ass{{\rm ass}}

\def\ker{ {\rm ker } }

\def\SL2Z{ {\rm SL}_2({\bf Z}) }

\def\th{ \theta }

\def\CL{{\cal L}}

\def\Gp1{ G^{1 , 1 } }
\def\P11{ P^{-1 , 1 } }
\def\Pp1{ P^{1 , 1 } }

\def\th{\theta}

\def\nCLsr{{}^\nu\kern-2pt {\cal L}^{\sigma , \rho  }}
\def\nP{{}^\nu \kern-2pt P}
\def\nL{{}^\nu\kern-2pt L}
\def\nLL{{}^\nu\kern-2pt \Lambda}
\def\nPsr{{}^\nu\kern-2pt P^{\sigma , \rho  }}
\def\nLsr{{}^\nu\kern-2pt L^{\sigma , \rho  }}
\def\nuCL{{}^\nu\kern-2pt  {\cal L}}
\def\nCLsr{{}^\nu\kern-2pt {\cal L}^{\sigma , \rho  }}
\def\nCL1m{{}^\nu\kern-2pt {\cal L}^{-1 , 1  }}
\def\x1nu{x^\frac{1}{\nu}}
\def\xm1nu{x^{-\frac{1}{\nu}}}

\def\CN{{\cal N}}
\def\ra{\rightarrow }

\def\CB{{\cal B}}

\def\CT{{\cal T}}

\def\CC{ {\cal C}}

\def\nAM0{{\nu }_{{\cal A},M_0}}
\def\nAN0{{\nu }_{{\cal A},N_0}}

\def\ga{\mathfrak{a}}
\def\gb{\mathfrak{b}}
\def\gc{\mathfrak{c}}

\def\SL{{\rm SL}}

\def\di!{\frac{\der^i}{i!}}
\def\dik!{\frac{\der^k_i}{k!}}

\def\N{\mathbb{N}}

\def\0{\overline{0}}
\def\1{\overline{1}}

\def\Ln1{\L_{n,\overline{1}}}

\def\a1{a_{\overline{1}}}

\def\S{\Sigma}

\def\vn1{\overrightarrow{n-1}}

\def\Min{{\rm Min}}

\def\mJ{\mathbb{J}}
\def\mI{\mathbb{I}}

\def\K1{{\rm K}_1}

\def\hmI1{\widehat{\mI_1}}
\def\tmI1{\widetilde{\mI_1}}
\def\tmJ1{\widetilde{\mJ_1}}
\def\hB1{\widehat{B_1}}
\def\hCB1{\widehat{\CB_1}}

\def\Den{{\rm Den}}
\def\Denl{{\rm Den}_l}
\def\Ore{{\rm Ore}}

\def\Den{{\rm Den}}

\def\Loc{{\rm Loc}}

\def\Ass{{\rm Ass}}

\def\maxDen{{\rm max.Den}}
\def\maxAss{{\rm max.Ass}}
\def\maxLoc{{\rm max.Loc}}
\def\llrad{{\rm l.lrad}}

\def\assmaxDen{{\rm ass.max.Den}}

\def\ga{\mathfrak{a}}

\def\udim{{\rm udim}}

\def\gll{\mathfrak{l}}

\begin{document}

\author{V. V. \  Bavula 
}

\title{Left localizable rings and their characterizations}

\maketitle

\begin{abstract}

A new class of rings, {\em the class of left localizable rings},  is
introduced. A ring $R$ is {\em left localizable} if each nonzero element
of $R$ is invertible in some left localization $S^{-1}R$ of the
ring $R$. Explicit criteria are given for  a ring to be a left
localizable ring provided the ring has only finitely many maximal
left denominator sets (eg, this is the case if a ring has a left
Artinian left quotient ring). It is proved that a ring with finitely
many maximal left denominator sets is a left localizable ring iff
 its left quotient
ring is a direct product of finitely many  division rings. A characterization is given of the class of rings that are finite direct product of left localization maximal rings.




$${\bf Contents}$$
\begin{enumerate}
\item Introduction.
\item Direct products of
left localization maximal rings and their characterization. \item
Left localizable rings and their characterization.
 \item The core of a left Ore set.
\end{enumerate}
\end{abstract}


\section{Introduction}


Throughout, module means a left module. In this paper the following notation will remained fixed.

$\noindent $

{\bf Notation}:

\begin{itemize}

\item $R$ is a ring with 1 and
$\CC_R$ is the set of (left and right) regular elements of the
ring $R$  (i.e. $\CC_R$ is the set of non-zero-divisors of $R$);
\item  $Q_{l, cl}(R):=\CC_R^{-1}R$ is the {\em classical left quotient ring} of $R$ (if it exists);

\item $\Ore_l(R):=\{ S\, | \, S$ is a left Ore set in $R\}$; \item
$\Den_l(R):=\{ S\, | \, S$ is a left denominator set in $R\}$;
\item $\Loc_l(R):= \{ [S^{-1}R]\, | \, S\in \Den_l(R)\}$ where $[S^{-1}R]$ is an $R$-{\em isomorphism} class of the ring $S^{-1}R$ (a ring isomorphism $\s : S^{-1}R\ra S'^{-1}R$ is called an $R$-{\em isomorphism} if $\s (\frac{r}{1})= \frac{r}{1}$ for all elements $r\in R$);
   \item $\gll_R:= \llrad (R)$ is the {\em left localization radical} of the ring $R$;
    \item
$\Ass_l(R):= \{ \ass (S)\, | \, S\in \Den_l(R)\}$ where $\ass
(S):= \{ r\in R \, | \, sr=0$ for some $s=s(r)\in S\}$; \item
$\Den_l(R, \ga ) := \{ S\in \Den_l(R)\, | \, \ass (S)=\ga \}$
where $\ga \in \Ass_l(R)$;
\item      $\Denl
(R,0)$ is the set of regular  left denominator sets $S$ in $R$
($S\subseteq \CC_R$);
\item $S_\ga=S_\ga (R)=S_{l,\ga }(R)$
is the {\em largest element} of the poset $(\Den_l(R, \ga ),
\subseteq )$ and $Q_\ga (R):=Q_{l,\ga }(R):=S_\ga^{-1} R$ is  the
{\em largest left quotient ring associated to} $\ga$, $S_\ga $
exists (Theorem 
 2.1, \cite{larglquot});
\item In particular, $S_0=S_0(R)=S_{l,0}(R)$ is the largest
element of the poset $(\Den_l(R, 0), \subseteq )$ and
$Q_l(R):=S_0^{-1}R$ is the largest left quotient ring of $R$;
\item $\Loc_l(R, \ga ):= \{ [S^{-1}R]\, | \, S\in \Den_l(R, \ga
)\}$.
\end{itemize}


{\bf Left localizable rings and their characterization}.  The aim of the paper is to introduce a new class of rings, the class of left localizable rings, and to give several  characterizations of them (Theorem \ref{29Nov12}, Theorem \ref{D2Dec12} and Theorem \ref{3Dec12}) in the case when they admit only finitely many maximal left denominator sets. Notice that each ring with left Artinian, left quotient ring has only finitely many maximal left denominator sets, \cite{Bav-LocArtRing}. A ring $R$ is called a {\em left localizable ring} if for  each nonzero element $r\in R$ there exists a left denominator set $S= S(r)$ such that $r\in S$ (equivalently, if for each nonzero element  $r\in R$ there exists a left denominator set $S'= S'(r)$ such that $\frac{r}{1}$ is a unit in $S'^{-1}R$).

\begin{itemize}
\item ({\bf Theorem \ref{29Nov12}})
 {\em  Let $R$ be a ring with $\maxDen_l(R)=\{ S_1, \ldots , S_n\}$. Let
$\ga_i:= \ass (S_i)$; $\s_i :R\ra R_i:=S_i^{-1}R$, $r\mapsto
\frac{r}{1}=r_i$; and $\s := \prod_{i=1}^n \s_i : R\ra
\prod_{i=1}^nR_i$, $r\mapsto (r_1, \ldots , r_n)$. The following
statements are equivalent.}
\begin{enumerate}
\item  {\em The ring $R$ is a left localizable ring.
\item $\gll_R=0$ and the rings $R_1, \ldots , R_n$ are division rings.
\item The homomorphism $\s$ is an injection and the rings $R_1, \ldots  , R_n$ are division rings. }
\end{enumerate}
\end{itemize}

\begin{itemize}
\item ({\bf Theorem \ref{3Dec12}}) {\em Let $R$ be a ring. The following statements are equivalent.}
\begin{enumerate}
\item {\em The ring $R$ is a left localizable ring with
$n:=|\maxDen_l(R)|<\infty$.\item $Q_{l,cl}(R)=R_1\times\cdots
\times R_n$  where $R_i$ are division rings. }
 \item {\em The ring $R$ is a semiprime left Goldie ring with $\udim
 (R)=|\Min (R)|=n$ where $\Min (R)$ is the set of minimal prime
 ideals of the ring $R$.}
\item $Q_l(R)= R_1\times\cdots \times R_n$ {\em  where $R_i$  are division rings.}
\end{enumerate}
\end{itemize}

{\bf Direct product of left localization maximal rings and their characterization}.  A ring $R$ is called a {\em left localization maximal ring} if every left denominator set of $R$ consists of units, i.e. one cannot invert anything new `on the left' in $R$, \cite{Crit-S-Simp-lQuot}. Theorem \ref{27Nov12} is a characterization of rings that are isomorphic to finite direct product of left localization maximal rings.

\begin{itemize}
\item ({\bf Theorem \ref{27Nov12}})
 {\em A ring $R$ is a direct product of $n$  left localization maximal rings iff the following conditions hold:}
\begin{enumerate}
\item  $\maxDen_l(R) = \{ S_1, \ldots , S_n\}$.
\item $\gll_R =0$.
\item $\ass (S_i) +\ass (S_j)=R$ {\em for all} $i\neq j$.
\item {\em The factor rings $R/\ass (S_i)$, $i=1, \ldots , n$, are  left localization maximal rings.}
\end{enumerate}
\end{itemize}

{\bf The core of a left Ore  set}. In Section \ref{TCLOS}, properties of the core of a left Ore set are studied. Let $R$ be a ring and $S$ be its left Ore set. The subset of $S$,
$$ S_c:=\{ s\in S\, | \, \ker (s\cdot ) = \ass (S)\}$$ is called the {\em core} of the left Ore set $S$ where $s\cdot : R\ra R$, $r\mapsto sr$, \cite{Crit-S-Simp-lQuot}.


\begin{itemize}
\item ({\bf Theorem \ref{A2Dec12}})
{\em Suppose that $S\in \Den_l(R, \ga )$ and $S_c\neq \emptyset$. Then}
\begin{enumerate}
\item $S_c\in \Den_l(R, \ga )$.
\item {\em The map $\th : S_c^{-1}R\ra S^{-1}R$, $s^{-1}r\mapsto s^{-1}r$, is a ring isomorphism. So,} $S_c^{-1}R\simeq S^{-1}R$.
\end{enumerate}
\end{itemize}

The next theorem gives an explicit description of the cores of maximal left denominator sets of a left localizable ring provided $|\maxDen_l(R)|<\infty$.
\begin{itemize}
\item ({\bf Theorem \ref{C2Dec12}})
{\em Let $R$ be a left localizable ring such that $\maxDen_l(R) = \{ S_1, \ldots , S_n\}$.}
\begin{enumerate}
\item {\em If $n=1$ then $S_{1,c} = S_1= R\backslash \{ 0\}$.
\item If $n\geq 2$ then $S_{i, c} = S_i\cap \bigcap_{j\neq i} \ga_i$  where $\ga_j=\ass (S_j)$.}
\end{enumerate}
\end{itemize}


\section{Direct products of left localization maximal rings and their characterization}\label{DPLLMC}

The aim of this section is to give a criterion of when a ring is isomorphic to a finite direct  product of left localization maximal rings (Theorem \ref{27Nov12}). At the beginning of the section necessary results are given  that are used in the proofs of the paper.

In \cite{larglquot}, we introduce the following new concepts and to
prove their existence for an {\em arbitrary} ring: {\em the
largest left quotient ring of a ring, the largest  regular left
 Ore  set of a ring,  the maximal left quotient rings of a ring,
the largest (two-sided) quotient ring of a ring, the maximal
(two-sided) quotient rings of a ring, a (left) localization
maximal ring}. Let us give more details.

{\bf The largest regular left Ore set and the largest left
quotient ring of a ring}. Let $R$ be a ring. A {\em
multiplicatively closed subset} $S$ of $R$ or a {\em
 multiplicative subset} of $R$ (i.e. a multiplicative sub-semigroup of $(R,
\cdot )$ such that $1\in S$ and $0\not\in S$) is said to be a {\em
left Ore set} if it satisfies the {\em left Ore condition}: for
each $r\in R$ and
 $s\in S$,
$$ Sr\bigcap Rs\neq \emptyset .$$
Let $\Ore_l(R)$ be the set of all left Ore sets of $R$.
  For  $S\in \Ore_l(R)$, $\ass (S) :=\{ r\in
R\, | \, sr=0 \;\; {\rm for\;  some}\;\; s\in S\}$  is an ideal of
the ring $R$.


A left Ore set $S$ is called a {\em left denominator set} of the
ring $R$ if $rs=0$ for some elements $ r\in R$ and $s\in S$ implies
$tr=0$ for some element $t\in S$, i.e. $r\in \ass (S)$. Let
$\Den_l(R)$ be the set of all left denominator sets of $R$. For
$S\in \Den_l(R)$, let $S^{-1}R=\{ s^{-1}r\, | \, s\in S, r\in R\}$
be the {\em left localization} of the ring $R$ at $S$ (the {\em
left quotient ring} of $R$ at $S$).

In general, the set $\CC$ of regular elements of a ring $R$ is
neither left nor right Ore set of the ring $R$ and as a
 result neither left nor right classical  quotient ring ($Q_{l,cl}(R):=\CC^{-1}R$ and
 $Q_{r,cl}(R):=R\CC^{-1}$) exists.
 Remarkably, there  exists the largest
 regular left Ore set $S_0= S_{l,0} = S_{l,0}(R)$. This means that the set $S_{l,0}(R)$ is an Ore set of
 the ring $R$ that consists
 of regular elements (i.e., $S_{l,0}(R)\subseteq \CC$) and contains all the left Ore sets in $R$ that consist of
 regular elements. Also, there exists the largest regular (left and right) Ore set $S_{l,r,0}(R)$ of the ring $R$.
 In general, all the sets $\CC$, $S_{l,0}(R)$, $S_{r,0}(R)$ and $S_{l,r,0}(R)$ are distinct, for example,
 when $R= \mI_1:=K\langle x, \frac{d}{dx}, \int\rangle$ is the ring of polynomial integro-differential operators over a field $K$ of characteristic zero,
  \cite{intdifline}.

$\noindent $

{\it Definition}, \cite{intdifline}, \cite{larglquot}.    The ring
$$Q_l(R):= S_{l,0}(R)^{-1}R$$ (respectively, $Q_r(R):=RS_{r,0}(R)^{-1}$ and
$Q(R):= S_{l,r,0}(R)^{-1}R\simeq RS_{l,r,0}(R)^{-1}$) is  called
the {\em largest left} (respectively, {\em right and two-sided})
{\em quotient ring} of the ring $R$.

$\noindent $

 In general, the rings $Q_l(R)$, $Q_r(R)$ and $Q(R)$
are not isomorphic, for example, when $R= \mI_1$, \cite{intdifline}.  The next
theorem gives various properties of the ring $Q_l(R)$. In
particular, it describes its group of units.


\begin{theorem}\label{4Jul10}
\cite{larglquot}
\begin{enumerate}
\item $ S_0 (Q_l(R))= Q_l(R)^*$ {\em and} $S_0(Q_l(R))\cap R=
S_0(R)$.
 \item $Q_l(R)^*= \langle S_0(R), S_0(R)^{-1}\rangle$, {\em i.e. the
 group of units of the ring $Q_l(R)$ is generated by the sets
 $S_0(R)$ and} $S_0(R)^{-1}:= \{ s^{-1} \, | \, s\in S_0(R)\}$.
 \item $Q_l(R)^* = \{ s^{-1}t\, | \, s,t\in S_0(R)\}$.
 \item $Q_l(Q_l(R))=Q_l(R)$.
\end{enumerate}
\end{theorem}

The set $(\Den_l(R), \subseteq )$ is a poset (partially ordered
set). In \cite{larglquot}, it is proved  that the set
$\maxDen_l(R)$ of its maximal elements is a {\em non-empty} set.

{\bf The maximal denominator sets and the maximal left localizations  of a ring}.


{\it Definition}, \cite{larglquot}. An element $S$ of the set
$\maxDen_l(R)$ is called a {\em maximal left denominator set} of
the ring $R$ and the ring $S^{-1}R$ is called a {\em maximal left
quotient ring} of the ring $R$ or a {\em maximal left localization
ring} of the ring $R$. The intersection
\begin{equation}\label{llradR}
\gll_R:=\llrad (R) := \bigcap_{S\in \maxDen_l(R)} \ass (S)
\end{equation}
is called the {\em left localization radical } of the ring $R$,
\cite{larglquot}.

 For a ring $R$, there is the canonical exact
sequence 
\begin{equation}\label{llRseq}
0\ra \gll_R \ra R\stackrel{\s }{\ra} \prod_{S\in \maxDen_l(R)}S^{-1}R, \;\; \s := \prod_{S\in \maxDen_l(R)}\, \s_S,
\end{equation}
where $\s_S:R\ra S^{-1}R$, $r\mapsto \frac{r}{1}$.

$\noindent $

{\bf The maximal elements of $\Ass_l(R)$}.  Let $\maxAss_l(R)$ be
the set of maximal elements of the poset $(\Ass_l(R), \subseteq )$
and

\begin{equation}\label{mADen}
\assmaxDen_l(R) := \{ \ass (S) \, | \, S\in \maxDen_l(R) \}.
\end{equation}
These two sets are equal (Proposition \ref{b27Nov12}), a proof is
based on Lemma \ref{1a27Nov12}. For a non-empty subset $X$ of $R$, let ${\rm r.ass} (X):=\{ r\in R\, | \, rx=0$ for some $x\in X\}$.

\begin{lemma}\label{1a27Nov12}
 \cite{larglquot}
Let $S\in \Den_l(R, \ga )$ and $T\in \Den_l(R, \gb )$ be such that $ \ga \subseteq \gb$. Let $ST$ be the multiplicative semigroup generated by $S$ and $T$ in $(R,\cdot )$.  Then
\begin{enumerate}
\item ${\rm r.ass} (ST)\subseteq \gb$.
\item $ST \in \Den_l(R, \gc )$ and $\gb \subseteq \gc$.
\end{enumerate}
\end{lemma}

\begin{corollary}\label{d4Jan13}
Let $R$ be a ring, $S\in \maxDen_l(R)$ and $T\in \Den_l(R)$. Then $T\subseteq S$ iff $\ass (T)\subseteq \ass (S)$.

\end{corollary}

{\it Proof}. $(\Rightarrow )$ If $T\subseteq S$ then $\ass (T)\subseteq \ass (S)$.

$(\Leftarrow )$ If $\ass (T)\subseteq \ass (S)$ then, by Lemma \ref{1a27Nov12}, $ST\in \Den_l(R)$ and $S\subseteq ST$, hence $S= ST$, by the maximality of $S$. Then $T\subseteq S$.  $\Box $


\begin{proposition}\label{b27Nov12}
\cite{larglquot} $\; \maxAss_l(R)= \assmaxDen_l(R)\neq \emptyset$. In particular, the ideals of this set are incomparable (i.e. neither $\ga\nsubseteq \gb$ nor $\ga\nsupseteq \gb$).
\end{proposition}

{\bf The localization maximal rings}. The set ($\Loc_l(R), \ra )$ is a partially ordered set (poset) where $A\ra B$ if there is an $R$-homomorphism $A\ra B$. There is no oriented loops in the poset $\Loc_l(R)$ apart from the $R$-isomorphism $A\ra A$.

Let $\maxLoc_l(R)$ be the set of maximal elements of the poset
$(\Loc_l(R), \ra )$. Then (see \cite{larglquot}),
\begin{equation}\label{mADen1}
\maxLoc_l(R) = \{ S^{-1}R \, | \, S\in \maxDen_l(R) \}= \{ Q_l(R/
\ga ) \, | \, \ga \in \assmaxDen_l(R)\}.
\end{equation}

 $\noindent $

 {\it Definition}, \cite{larglquot}. A ring $A$ is
called a {\em left localization maximal ring} if $A= Q_l(A)$ and
$\Ass_l(A) = \{ 0\}$. A ring $A$ is called a {\em right
localization maximal ring} if $A= Q_r(A)$ and $\Ass_r(A) = \{
0\}$. A ring $A$ which is a left and right localization maximal
ring is called a {\em (left and right) localization maximal ring}
(i.e. $Q_l(A) =A=Q_r(A)$ and $\Ass_l(A) =\Ass_r(A) = \{ 0\}$).

$\noindent $

The next theorem is a criterion of  when a left quotient ring of a
ring is a maximal left quotient ring of the ring.

\begin{theorem}\label{21Nov10}
\cite{larglquot} Let  a ring $A$ be a left localization of a ring
$R$, i.e. $A\in \Loc_l(R, \ga )$ for some $\ga \in \Ass_l( R)$.
Then $A\in \maxLoc_l( R)$ iff $Q_l( A) = A$ and  $\Ass_l(A) = \{
0\}$, i.e. $A$ is a left localization maximal ring (clearly, $\ga \in \assmaxDen_l(R)$).
\end{theorem}


Theorem \ref{21Nov10} shows that the left localization maximal
rings are precisely the localizations of all the rings at their
maximal left denominators sets.


{\it Example}. Let $A$ be a simple ring. Then $Q_l(A)$ is a left
localization maximal  ring and $Q_r(A)$ is a right localization
maximal ring (by Theorem \ref{4Jul10}.(4) and Theorem \ref{21Nov10}).


{\it Example}. A division ring is a (left and right) localization
maximal ring. More generally, a simple Artinian ring (i.e. the
matrix algebra over a division ring) is a (left and right)
localization maximal ring.


{\bf Left (non-)localizable elements of a ring}.

\begin{lemma}\label{b11Dec12}
\cite{Bav-LocArtRing} Let $S\in \Den_l(R, \ga )$ (respectively, $S\in \Den (R, \ga )$),
$\s : R\ra S^{-1}R$, $r\mapsto \frac{r}{1}$, and $G:=(S^{-1}R)^*$
be the group of units of the ring $S^{-1}R$. Then $S':=\s^{-1}
(G)\in \Den_l(R, \ga )$ (respectively, $S':=\s^{-1} (G)\in \Den
(R, \ga )$).
\end{lemma}


{\it Definition}, \cite{Crit-S-Simp-lQuot}. An element $r$ of a ring $R$ is called a {\em left localizable element} if there  exists a left denominator set $S$  of
 $R$ such that $r\in S$ (and so the element $\frac{r}{1}\neq 0$ is invertible in the ring
 $S^{-1}R$), equivalently,  if there  exists a left denominator set $T$  of
 $R$ such that the element $\frac{r}{1}$ is invertible in the ring
 $T^{-1}R$ (Lemma \ref{b11Dec12}). The set of left localizable elements is denoted $\CL_l(R)$.

$\noindent $

Clearly, 
\begin{equation}\label{CLeU}
\CL_l(R)=\bigcup_{S\in \maxDen_l(R)} S.
\end{equation}

Similarly, a {\em right localizable element} is defined and let
$\CL_r(R)$ be the set of right localizable elements of the ring
$R$. The elements of the {\em set of left and right localizable
elements},
$$ \CL_{l,r}(R) = \CL_l(R)\cap \CL_r(R),$$
 are called {\em left and right localizable elements}. An element $r\in R$ is called a {\em localizable element} if there exists a
  (left and right) denominator set $S\in \Den (R)$ such that $r\in S$, equivalently,  if there exists a
  (left and right) denominator set $T\in \Den (R)$ such that  the element $\frac{r}{1}$ is invertible in the ring
 $T^{-1}R$ (Lemma \ref{b11Dec12}).  The set of all localizable elements of the ring $R$ is denoted by $\CL (R)$. Clearly,
 $$ \CL (R) \subseteq \CL_{l,r}(R).$$
 The sets
 $$ \CN\CL_l(R) := R\backslash \CL_l (R), \;\; \CN\CL_r(R) := R\backslash \CL_r (R), \;\;
 \CN\CL_{l,r}(R) := R\backslash \CL_{l,r} (R), \;\;\CN\CL (R) := R\backslash \CL (R),$$

are called the {\em  sets of left non-localizable; of right
non-localizable; of left and right non-localizable; of
non-localizable elements}, respectively. The elements of these
sets are called correspondingly (eg, an element $r\in \CN\CL_l(R)$
is called a {\em left non-localizable element}).

\begin{lemma}\label{a27Nov12}

\begin{enumerate}
\item $\gll_R\cap \CL_l (R)=\emptyset$.
\item $\gll_R \subseteq \CN \CL_l(R)$.
\end{enumerate}
\end{lemma}

{\it Proof}. 1. If $s\in \CL_l(R)$ then $s\in S$ for some $S\in
\maxDen_l(R)$ (by (\ref{CLeU})), and so $\frac{s}{1}\neq 0$ in
$S^{-1}R$. Therefore, $s\not\in \gll_R$.

2. Statement 2 follows from statement 1 as $\CN\CL_l(R)= R\backslash \CL_l(R)$. $\Box $

$\noindent $

{\it Definition}, \cite{Crit-S-Simp-lQuot}. For an arbitrary ring $R$, the intersection $$\CC_l(R):=\bigcap_{S\in \maxDen_l(R)}S$$
is called the set of {\em completely left localizable elements} of $R$ and an element of the set $\CC_l(R)$ is called a {\em completely left  localizable element}.

$\noindent $

{\bf The maximal left quotient rings of a finite direct product of rings}.
\begin{theorem}\label{c26Dec12}
\cite{Crit-S-Simp-lQuot}  Let $R=\prod_{i=1}^n R_i$ be a direct product of rings $R_i$. Then
for each $i=1, \ldots , n$, the map
\begin{equation}\label{1aab1}
\maxDen_l(R_i) \ra \maxDen_l(R), \;\; S_i\mapsto R_1\times\cdots \times S_i\times\cdots \times R_n,
\end{equation}
is an injection. Moreover, $\maxDen_l(R)=\coprod_{i=1}^n \maxDen_l(R_i)$ in the sense of (\ref{1aab1}), i.e.
$$ \maxDen_l(R)=\{ S_i\, | \, S_i\in \maxDen_l(R_i), \; i=1, \ldots , n\},$$
$S_i^{-1}R\simeq S_i^{-1}R_i$, $\ass_R(S_i)= R_1\times \cdots \times \ass_{R_i}(S_i)\times\cdots \times R_n$. The core of the left denominator set $S_i$ in $R$ coincides with the core $S_{i,c}$ of the left denominator set $S_i$ in $R_i$, i.e.
$$(R_1\times\cdots \times S_i\times\cdots \times R_n)_c=0\times\cdots \times S_{i,c}\times\cdots \times 0.$$
\end{theorem}


{\bf Finite direct products of localization maximal rings}.  The
next theorem describes properties of finite direct products of
left localization maximal rings.  It is used in the proof of
Theorem \ref{27Nov12} which is an explicit characterization of all
such rings.

\begin{theorem}\label{25Nov12}
Let $R_1, \ldots , R_n$ be left localization maximal rings and $R:= R_1\times \cdots \times R_n$. Then
\begin{enumerate}
\item $\maxDen_l(R) = \{ S_i\, | \, i=1, \ldots , n\}$ and $S_i:= R_1\times \cdots \times R_i^*\times \cdots \times R_n$ where $R_i^*$ is a the group of units of the ring $R_i$.
\item $\ass (S_i)=R_1\times \cdots \times 0 \times \cdots \times R_n$, $i=1, \ldots , n$ where $0$ is on $i$'th place; $\ass (S_i) +\ass (S_j) = R$ for all $i\neq j$.
\item $S_i^{-1}R\simeq R_i\simeq R/\ass (S_i)$, $i=1, \ldots , n$. The map $\s_i : R\ra S_i^{-1}R$, $ r\mapsto \frac{r}{1}$, is equal to the map $R\ra R/ \ass (S_i)$, $ r\mapsto r+\ass (S_i)$.
    \item $\gll_R =0$.
    \item The natural ring homomorphism
    $$ \s := \prod_{I=1}^n \s_i: R\ra \prod_{I=1}^n S_i^{-1}R, \;\; r\mapsto ( \frac{r}{1}, \ldots ,  \frac{r}{1})=(r_1, \ldots , r_n), $$
    is an isomorphism.
    \item $\CC_l(R)=\bigcap_{i=1}^n S_i = R^*$ where $R^*= R_1^*\times \cdots \times R_n^*$ is the group of units of the ring $R$.
        \item The set $\CL_l(R)=\bigcup_{i=1}^n S_i$ of left localizable elements of the ring $R$ contains precisely the elements $(r_1, \ldots , r_n)\in R$ such that $r_i\in R_i^*$ for some $i$.
            \item The set $\CN\CL_l(R)$ of left non-localizable elements of the ring $R$ is equal to $\prod_{i=1}^n R_i\backslash R_i^*$, the direct product of the sets of left non-localizable elements $\CN\CL_l(R_i)=R_i\backslash R_i^*$ of the rings $R_i$.
\end{enumerate}
\end{theorem}

{\it Proof}. 1. Statement 1 is a particular case of Theorem \ref{c26Dec12}.

 2-5. For $i=1, \ldots , n$,
$$ S_i\in \maxDen_l(R), \;\; \ass (S_i) = R_1 \times \cdots  \times 0 \times \cdots \times R_n , \;\; S_i^{-1}R\simeq R_i \simeq R/\ass (S_i).$$
and the localization map $\s_i$ is equal to the map $R\ra R/ \ass (S_i)$, $ r\mapsto r+\ass (S_i)$. Clearly, $\ass (S_i) +\ass (S_j) = R$ for all $i\neq j$ and $\gll_R \subseteq \cap_{i=1}^n \ass (S_i) =0$, and so $\gll_R =0$. The ring homomorphism $\s $ is an isomorphism. So, statements 2-5 have been proved.

6. Statement 6 follows from statement 1.

7. Trivial.

8. Statement 8 follows from statement 1. $\Box $


\begin{corollary}\label{a29Nov12}
Let $R=\prod_{i=1}^n R_i$ be a direct product of left localization maximal rings $R_i$. Then every nonzero element of the ring $R$ is left localizable iff the ring $R_1, \ldots , R_n$ are division rings.
\end{corollary}

{\it Proof}. Every nonzero element of the ring $R$ is left localizable iff $\{ 0\} = \CN\CL_l(R) = \prod_{i=1}^n R_i\backslash R_i^*$ (Theorem \ref{25Nov12}.(8)) iff $R_i\backslash R_i^*=\{0\}$ for all $i$ iff the ring $R_1, \ldots , R_n$ are division rings. $\Box $

$\noindent $

Corollary \ref{a29Nov12} is a particular case of a criterion
(Theorem \ref{29Nov12}) for a ring with finitely many maximal left
denominator sets to be a left localizable ring.


{\bf Characterization of finite products of left localization
maximal rings}. The next theorem is a criterion for a ring to be
isomorphic to a finite direct  product of left localization
maximal rings.

\begin{theorem}\label{27Nov12}
 A ring $R$ is a direct product of $n$  left localization maximal rings iff the following conditions hold:
\begin{enumerate}
\item  $\maxDen_l(R) = \{ S_1, \ldots , S_n\}$.
\item $\gll_R=0$.
\item $\ass (S_i) +\ass (S_j)=R$ for all $i\neq j$.
\item The factor rings $R/\ass (S_i)$, $i=1, \ldots , n$, are  left localization maximal rings.
\end{enumerate}
In this case,

(a)  for each $i=1, \ldots , n$, the ring $R_i:= S_i^{-1}R$ is a
left localization maximal ring such that the homomorphism $\s_i :
R\ra R_i$, $ r\mapsto \frac{r}{1}$, is an epimorphism with $\ker
(\s_i) = \ass (S_i)$, and so $R_i\simeq R/\ass (S_i)$.

(b)  The
map $\s : =\prod_{i=1}^n \s_i : \; R\ra \prod_{i=1}^n R_i$ is an
isomorphism.

(c)  Let us identify the ring $R$ with
$\prod_{i=1}^n R_i$ via $\s$. Then $S_i = R_1\times \cdots \times
R_i^* \times \cdots \cdots \times R_n=\s_i^{-1}(R_i^*)$ where
$R_i^* = (R/ \ass (S_i))^*$ is the group of units of the ring
$R_i$, and for each $i=1, \ldots , n$, $R_i^* = \s_i (S_i)$.

(d) $\ass (S_i) = R_1\times \cdots \times 0\times \cdots \times R_n$ for $i=1,\ldots , n$ where $0$ is on $i$'th place.

\end{theorem}

{\it Proof}.  $(\Rightarrow )$ Theorem \ref{25Nov12}.

$(\Leftarrow )$  By (\ref{llRseq}) and conditions 1 and 2,  the
map $\s $ is a ring monomorphism since
$$ \ker (\s ) = \bigcap_{i=1}^n \ass (S_i) = \gll_R =0.$$
 The ring homomorphism $\s_i$ yields a monomorphism $R/ \ass (S_i) \ra S_i^{-1}R= R_i$ which is necessarily an isomorphism, by condition 4. Clearly, $\s_i: R\ra R/ \ass (S_i)$, $r\mapsto r+\ass (S_i)$.

Conditions 1 and 2, $\bigcap_{i=1}^n \ass (S_i) =0$ and $\ass
(S_i) +\ass (S_j) = R$ for all $i\ne j$, imply that the map $\s$
is an isomorphism, i.e. the ring $R\simeq R_1\times \cdots \times
R_n$ is a direct product of $n$ left localization maximal rings.

By Theorem \ref{25Nov12}.(1), $S_i=R_1\times \cdots \times
R_i^*\times \cdots \times R_n=\s_i^{-1}(R_i^*)$ and $\ass (S_i) =
R_1\times \cdots \times  0\times \cdots \times R_n$. Then $R_i^* =
\s_i (S_i)$ since $\s_i : R\ra R/\ass (S_i)$, $r\mapsto r+\ass
(S_i)$. $\Box $

$\noindent $


\section{Left localizable rings and their characterization}\label{LLRTC}

The aim of this section is to introduce a new class of rings, {\em
the class of left localizable rings}, consider their properties,
and give their  characterizations   (Theorem \ref{29Nov12}, Theorem \ref{D2Dec12} and
Theorem \ref{3Dec12}) provided they have only finitely many
maximal left denominator sets. In particular, we prove that  a
ring with finitely many maximal left denominator sets is a left
localizable ring iff it is a semiprime left Goldie ring such that
its left quotient ring is a direct product of division rings
(Theorem \ref{3Dec12}). For a semiprime left Goldie ring $R$, it
is proved that $\maxDen_l(R)|<\infty$. Moreover, elements of
$\maxDen_l(R)$ are found (Corollary \ref{a4Dec12}).

$\noindent $

{\it Definition}.  A ring $R$ is called a {\em left localizable ring} if the set of left localizable elements $\CL_l(R)$ is
 the largest possible, i.e. $\CL_l(R)=R\backslash \{ 0\}$; equivalently, for every nonzero element $r\in R$
 there is a left denominator set $S= S(r)\in \Den_l(R)$ such that $r\in S$; equivalently, for every
 nonzero element $r\in R$ there is a left denominator set $S'= S'(r)$ such that the element $\frac{r}{1} \in S'^{-1}R$ is a
  unit of the ring $S'^{-1}R$ (Lemma \ref{b11Dec12}).
$\noindent $

Similarly,   right localizable; left and right localizable; and localizable rings are defined, i.e. where the sets $\CL_r(R)$; $\CL_{l,r}(R)$;  and $\CL(R)$  are equal to $R\backslash \{ 0\}$,  respectively. A localizable ring is also called a {\em two-sided localizable ring}.

\begin{theorem}\label{8Feb13}
Let $R=\prod_{i=1}^n R_i$ be a direct product of rings $R_i$. The ring $R$ is a left localizable ring iff the rings $R_i$ are so.
\end{theorem}

{\it Proof}. The proof is an easy corollary of Theorem \ref{c26Dec12} that states that $\maxDen_l(\prod_{i=1}^n R_i)= \coprod_{i=1}^n \maxDen_l(R_i)$.

$(\Rightarrow )$ Suppose that the ring $R$ is a left localizable ring. We have to show that the rings $R_i$ are so. Each  ring $R_i$ is a subring of $R$. Let $r_i\in R_i$ be a nonzero element. The ring $R$ is a left localizable ring. So, $r_i\in S_i$ for some $S_i\in \maxDen_l(R)$. Then $S_i\in \maxDen_l(R_i)$, by Theorem \ref{c26Dec12}.

$(\Leftarrow )$ Suppose that the rings $R_i$ are left localizable rings. Let $r= (r_1, \ldots , r_n)\in R$ be a nonzero element. Then $0\neq r_i\in R_i$ for some $i$. The ring $R_i$ is a left localizable ring, and so $r_i\in S_i$ for some $S_i\in \maxDen_l(R_i)$. By Theorem \ref{c26Dec12}, $S_i\in \maxDen_l(R)$. Therefore, $R$ is a left localizable ring. $\Box $

$\noindent $

\begin{lemma}\label{a1Dec12}
\begin{enumerate}
\item If a ring $R$ is a left (respectively, right; left and
right; two-sided) localizable then $\gll_R=0$ (respectively,
${\rm r.lrad}(R)=0$;  ${\rm lr.lrad}(R)=0$; ${\rm lrad}(R)=0$).
\item Let a ring $R= \prod_{i=1}^n R_i$ be a direct product of
left localization maximal rings $R_i$. Then $R$ is a left
localizable ring iff the rings $R_i$ are division rings.
\end{enumerate}
\end{lemma}

{\it Proof}. 1. Statements are easy corollaries of Lemma
\ref{a27Nov12}.

2. Corollary \ref{a29Nov12} $\Box $


{\bf Properties of the maximal left quotient rings of a ring}.
The next theorem describes various properties of the maximal left
quotient rings of a ring, in particular, their groups of units and
their largest left quotient rings.

\begin{theorem}\label{15Nov10}
\cite{larglquot} Let $S\in \maxDen_l(R)$, $A= S^{-1}R$, $A^*$ be
the group of units of the ring $A$; $\ga := \ass (S)$, $\pi_\ga
:R\ra R/ \ga $, $ a\mapsto a+\ga$, and $\s_\ga : R\ra A$, $
r\mapsto \frac{r}{1}$. Then
\begin{enumerate}
\item $S=S_\ga (R)$, $S= \pi_\ga^{-1} (S_0(R/\ga ))$, $ \pi_\ga
(S) = S_0(R/ \ga )$ and $A= S_0( R/\ga )^{-1} R/ \ga = Q_l(R/ \ga
)$. \item  $S_0(A) = A^*$ and $S_0(A) \cap (R/ \ga )= S_0( R/ \ga
)$. \item $S= \s_\ga^{-1}(A^*)$. \item $A^* = \langle \pi_\ga (S)
, \pi_\ga (S)^{-1} \rangle$, i.e. the group of units of the ring
$A$ is generated by the sets $\pi_\ga (S)$ and $\pi_\ga^{-1}(S):=
\{ \pi_\ga (s)^{-1} \, | \, s\in S\}$. \item $A^* = \{ \pi_\ga
(s)^{-1}\pi_\ga ( t) \, |\, s, t\in S\}$. \item $Q_l(A) = A$ and
$\Ass_l(A) = \{ 0\}$.     In particular, if $T\in \Den_l(A, 0)$
then  $T\subseteq A^*$.
\end{enumerate}
\end{theorem}


\begin{lemma}\label{d1Dec12}
Let $R$ be a ring and $S\in \maxDen_l(R)$. Then the ring $S^{-1}R$ is a division ring iff $R= S\cup \ass (S)$.
\end{lemma}

{\it Remark}. For any left Ore set $T$, $T\cap \ass (T) =\emptyset$.

{\it Proof}. Let $\s : R\ra S^{-1}R$, $ r\mapsto\frac{r}{1}$, and
$(S^{-1}R)^*$ be the group of units of the ring $S^{-1}R$. By
Theorem \ref{15Nov10}.(3), $ S= \s^{-1} ((S^{-1}R)^*)$  for all
$S\in \maxDen_l(R)$.

$(\Rightarrow )$ Let $s\in R$. Suppose that $s\not\in \ass (S)$, we have to show that $s\in S$. Clearly, $\frac{s}{1}\neq 0$ in $S^{-1}R$. Therefore,  $\frac{s}{1}\in (S^{-1}R)^*$ since $S^{-1}R$ is a division ring. Then  $s\in \s^{-1} ((S^{-1}R)^*)=S$.

$(\Leftarrow )$  If $R=S\cup \ass (S)$ then this union is a
disjoint union since $S\cap \ass (S) =\emptyset$. In particular,
$S+\ass (S)\subseteq S$. Let $\frac{r}{1}$ (where $r\in R$) be a
nonzero element of the ring $S^{-1}R$. Then $ r\not\in \ass (S)$,
and so  $r\in S$, i.e. $\frac{r}{1}$ is a unit in the ring
$S^{-1}R$. Therefore, the ring $S^{-1}R$ is a division ring. $\Box
$


{\bf Characterization of left localizable rings}. The next theorem
gives a characterization of left localizable rings with finite
number of maximal left denominator sets.

\begin{theorem}\label{29Nov12}
Let $R$ be a ring with $\maxDen_l(R)=\{ S_1, \ldots , S_n\}$. Let
$\ga_i:= \ass (S_i)$; $\s_i :R\ra R_i:=S_i^{-1}R$, $r\mapsto
\frac{r}{1}=r_i$; and $\s := \prod_{i=1}^n \s_i : R\ra
\prod_{i=1}^nR_i$, $r\mapsto (r_1, \ldots , r_n)$. The following
statements are equivalent.
\begin{enumerate}
\item The ring $R$ is a left localizable ring.
\item $\gll_R=0$ and the rings $R_1, \ldots , R_n$ are division rings.
\item The homomorphism $\s$ is an injection and the rings $R_1, \ldots  , R_n$ are division rings.
\end{enumerate}
\end{theorem}

{\it Proof}.  The rings $R_i$ are left localization maximal
(Theorem \ref{15Nov10}.(6)).

  $(2\Leftrightarrow 3)$ Statements 2 and 3 are equivalent since $\ker (\s ) = \bigcap_{i=1}^n \ass (S_i) = \gll_R$.

$(1\Leftarrow 2)$ Since $\gll_R=0$, the map $\s$ (see (\ref{llRseq})) is a ring monomorphism. So, every element $r$ is
a unique $n$-tuple $(r_1, \ldots , r_n)$ where $r_i\in R_i$. The rings $R_1, \ldots , R_n$ are division rings.
Let $r=(r_i)\in R$ be a nonzero element of $R$, then $r_i\neq 0$ for some $i$, and so $\s_i(r) = r_i\in R_i^*$.
Therefore, every  nonzero element of the ring is left localizable.

$(1\Rightarrow 2)$ Suppose that every nonzero element of the ring
$R$ is left localizable, i.e. $\CN\CL_l(R)=\{ 0\}$. By Lemma
\ref{a27Nov12}.(2), $\gll_R =0$, and so, by (\ref{llRseq}),
the map $\s$ is a ring monomorphism. We identify the ring $R$ with
its image $\s (R)$ in the ring $\prod_{i=1}^n R_i$, i.e. $r=(r_1,
\ldots , r_n)$ where $r_i= \s_i(r)$.

If $n=1$ then $\s_1: R\ra S_1^{_-1}R$, $r\mapsto \frac{r}{1}$, is
a monomorphism, and so $\ass (S_1)=0$. By the assumption
$S_1=\CL_l(R)=R\backslash \{ 0\}$, i.e.  $S_1^{-1}R$ is a division
ring.

Let $n\geq 2$. By Proposition \ref{b27Nov12}, the ideals $\ga_i :=
\ass (S_i)$, $i=1, \ldots , n$, are non-comparable, i.e.
\begin{equation}\label{aiajf}
\ga_i\not\subseteq \ga_j, \;\; i\neq j.
\end{equation}
In particular, $\ga_i \neq 0$ for all $i$. The proof consists of several steps.

{\em Step 1: For all $i$},
$$S_i\cap \bigcap_{j\neq i} \ga_j\neq  \emptyset .$$
Fix $r\in S_i$. In
particular, $r= (r_1, \ldots , r_n)$ with $r_i\in R_i^*$, and vice
versa, i.e. if $r_i\in R_i^*$ then $r\in \s_i^{-1} (R_i^*)=S_i$.
Fix a nonzero element, say $r\in S_i$, such that the $n$-tuple
$r=(r_1, \ldots , r_n)$ has the largest possible number of zeros,
say $m=n-s$ for some $s$. We claim that $s=1$. Up to order, we may
assume that $i=1$,
$$ r= (r_1, \ldots , r_s, 0, \ldots , 0)$$
and  $r_j\neq 0$ for $j=1, \ldots , s$, and $r_1\in R_1^*$. Suppose that $s\neq1$, we seek a contradiction. Then, for all elements $a_s\in \ga_s$,
$$\s (a_s) \s (r) = (\s_1(a_s)r_1, \ldots ,\s_{s-1} (a_s) r_{s-1}, 0, \ldots , 0).$$
By the maximality of $m$ we must have $\s(a_s) \s (r) =0$. In particular, $\s_1(a_s) r_1=0$, and so $\s_1(a_s)=0$ since $r_1\in R^*$, i.e. $\ga_s\subseteq \ga_1$, a contradiction,  (Proposition \ref{b27Nov12}).  Therefore, $s=1$, and $\s (r) \in S_1\cap \bigcap_{j\neq 1} \ga_j\neq 0$. The proof of Step 1 is complete.

 {\em Step 2: For each  $i$,  the ring $R_i$ is a division ring.}

 Without loss of generality we may assume that $i=1$, then   $S_1\cap \bigcap_{j=2}^n \ga_j\neq \emptyset$ (Step 1).
  Fix an element $s\in S_1\cap \bigcap_{j=2}^n \ga_j$. Then $\s (s) = (s_1, 0, \ldots , 0)$ and $s_1\in R_1^*$ (since $s\in S_1$ and $R_1=S_1^{-1}R$).
  Suppose that the ring $R_1$ is not a division ring, we seek a contradiction. Then, by Lemma \ref{d1Dec12}, $R_1\neq S_1\cup \ga_1$.
   Fix an element $r\in R\backslash (S_1\cup \ga_1)$. Then $\s (r) = (r_1, \ldots , r_n)$ with $ S_1^{-1}R\ni r_1\neq 0$ (since $r\not\in \ga_1$).
   Then
 $$ \s (sr) = (s_1r_1, 0, \ldots , 0) \in \prod_{i=1}^nR_i$$ with $s_1r_1\neq 0$ since $s_1\in R_1^*$ and $r_1\neq 0$.
 Therefore, $sr\not\in \bigcup_{j=2}^n S_j$ (since $sr\in \cap_{j=2}^n\ga_j$, by the choice of $s$). Then necessarily $sr\in S_1$ since $sr \in R\backslash \{ 0\}= \bigcup_{i=1}^n S_i$
 (the ring $R$ is a left localizable ring). Then $\s_1(r)= \s_1(s)^{-1}\s_1(sr) \in R_1^*$, and so $r\in \s_1^{-1}(R_1^*) = S_1$,
  by Theorem \ref{15Nov10}.(3). This contradicts to the fact that $r\in R\backslash (S_1\cup \ga_1)$. Therefore, $R_1$ is a division ring. The proof of Step 2 and the theorem  is complete.  $\Box $

\begin{corollary}\label{e1Dec12}
Let $R$ be a left localizable ring  with $\maxDen_l(R)=\{ S_1, \ldots , S_n\}$ where $n\geq 2$. Then, for all $i=1, \ldots , n$,
$S_i\cap \bigcap_{j\neq i}\ga_j= (\bigcap_{j\neq i}\ga_j)\backslash \{ 0\} \neq \emptyset $ where $\ga_i :=\ass (S_i)$.
\end{corollary}

{\it Proof}. By Step 1 of the proof of Theorem \ref{29Nov12},
 $S_i\cap \bigcap_{j\neq i}\ga_j\neq \emptyset $.  The equality in the corollary  follows from the
 fact that $\bigcap_{i=1}^n \ga_j=0$ and $R_i= S_i \cup \ga_i$ is a disjoint union (Lemma \ref{d1Dec12}, since $R_i$ is a division ring, by Theorem \ref{29Nov12}). Then   $(\bigcap_{j\neq i}\ga_j)\backslash \{ 0\} \neq \emptyset $. $\Box $


\begin{proposition}\label{A3Dec12}
Let $R$ be a left localizable ring with  $\maxDen_l(R)=\{ S_1,
\ldots , S_n\}$ where $n\geq 2$. We keep the notation of Theorem
\ref{29Nov12}. Let $C_i':=(\bigcap_{j\neq i}\ga_j)\backslash \{
0\}$. Then
\begin{enumerate}
\item $C':=C_1'+\cdots + C_n'\in \Den_l(R, 0)$ and $C'^{-1}R\simeq
R_1\times\cdots \times R_n$. \item For $i=1, \ldots, n$, $C_i'\in
\Den_l(R, \ga_i)$ and $C_i'^{-1}R=R_i$.
\end{enumerate}

\end{proposition}

{\it Proof}. Clearly, all $C_i'\neq \emptyset$ (Corollary
\ref{e1Dec12}); $C_i'C_j'=0$ for all $i\neq j$ (since
$C_i'C_j'\subseteq \bigcap_{k=1}^n \ga_k =0$, Theorem
\ref{29Nov12}), and $\s (C_i') \subseteq 0\times \cdots \times  0
\times R_i^*\times 0 \times\cdots \times 0$ since $C_i'=S_i\cap
\bigcap_{j\neq i} \ga_j$ (Corollary \ref{e1Dec12}). Therefore, the
set $C'$ consists of regular elements of the ring $R$ (since the
map $\s : R\ra \prod_{i=1}^n R_i$ is a monomorphism) and is
obviously multiplicatively closed. Let us identify the set $C_i'$
with its image in the division ring $R_i$ via $\s $. Then $\ass
(C_i') = \ga_i$ (since $\CC_i'S_i\subseteq \CC_i'$).  Notice that $C_i'\cup \{ 0\} =\bigcap_{j\neq i}
\ga_j$ is an {\em ideal} of the ring $R$. Fix an element $c_i\in
C_i'$, then each element $s^{-1}r\in R_i$, where $s\in S_i$ and
$r\in R$, can be written as a left fraction
$$ s^{-1} r= s^{-1} c_i^{-1} c_i r= (c_is)^{-1} c_i r, \;\; {\rm where}\;\; c_is\in C_i' \;\; {\rm and}\;\; c_ir \in C_i'\cup \{ 0\}.$$
Therefore,
 $C_i'\in \Den_l(R, \ga_i)$ and $ C_i'^{-1}R=R_i$, and statement 2 holds. So, any element $\alpha$ of the ring $\prod_{i=1}^n R_i$ can be written as $\alpha = (c_1^{-1}r_1, \ldots c_n^{-1}r_n)$
where $c_i\in C_i'$ and $r_i\in C_i'\cup \{ 0\}$, and so
$$\alpha = (c_1^{-1}r_1, \ldots c_n^{-1}r_n)=(c_1+\cdots + c_n)^{-1}(r_1+\cdots + r_n)$$ where
 $c:= c_1+\cdots + c_n\in C'$ and $r:= r_1+\cdots + r_n\in (\sum_{i=1}^n C_i') \cup \{ 0\} \subseteq R$. So,
 every element $\alpha$ of the ring $\prod_{i=1}^n R_i$ can be written as $\alpha = c^{-1} r$ for some elements
  $c\in C'$ and $r\in R$. Since the set $C'$ consists of regular elements, the fact above ($\alpha = c^{-1} r$)
  is equivalent to the fact that $C'$ is a  left Ore set, and so $C'\in \Den_l(R, 0)$. Clearly, $C'^{-1}R= R_1\times \cdots \times R_n$.  $\Box $

$\noindent $

The sets $C_i'$ in Theorem \ref{A3Dec12} have remarkable
properties: $C_i'\subseteq S_i$, $C_i'^{-1}R=S_i^{-1}R$ and $\ass
(c_i'\cdot ) =\ass (S_i)$ for all elements $ c_i'\in C_i'$ where
$c_i'\cdot : R\ra R$, $ r\mapsto c_i'r$ (Theorem \ref{C2Dec12}). In fact, the set $C_i'$ is
the core of $S_i$ (Theorem \ref{C2Dec12}).

The next theorem is a useful criterion for  a ring $R$ to be a
left localizable ring with $|\maxDen_l(R)|<\infty$.

\begin{theorem}\label{D2Dec12}
Let $R$ be a ring. The following statements are equivalent.
\begin{enumerate}
\item The ring $R$ is a left localizable ring with
$|\maxDen_l(R)|<\infty$. \item There are left denominator sets
$S_1', \ldots ,S_n'\in \Den_l(R)$ such that the rings $R_i:=
S_i'^{-1}R$, $i=1, \ldots , n$, are division rings and the map
    $$ \s : =\prod_{i=1}^n \s_i : R\ra \prod_{i=1}^n R_i, \;\; r\mapsto (\s_1(r), \ldots , \s_n(r)), $$ is an injection  where $\s_i : R\ra R_i$, $ r\mapsto \frac{r}{1}$.
\end{enumerate}
If one of the equivalent conditions holds and none of the rings
$R_i$ can be dropped with    preserving the injectivity of $\s$ then
$n = |\maxDen_l(R)|$, and $\maxDen_l(R)=\{ S_1, \ldots , S_n\}$
where $S_i= \s_i^{-1}(R_i^*)$.
\end{theorem}

{\it Proof}. $(1\Rightarrow 2)$ Theorem \ref{29Nov12}.

$(2\Leftarrow 1)$  Suppose that statement 2 holds and we assume
that none of the rings $R_i$ can be dropped, i.e. for each $i=1,
\ldots , n$, $\bigcap_{j\neq i} \ga_i \neq 0$ where $\ga_i := \ass
(S_i')$. We are going to prove that the implication $(2\Leftarrow
1)$ holds and that $\maxDen_l(R) = \{ S_1, \ldots , S_n\}$.

Every division ring is a left localization maximal ring. By
(Theorem 3.12.(3), \cite{larglquot}), $S_i\in \maxDen_l(R)$.
Clearly, $S_i'\subseteq \s_i^{-1} (R_i^*)= S_i$ and the natural
ring homomorphism $R_i=S_i'^{-1}R\ra S_i^{-1}R$,
$\frac{r}{1}\ra\frac{r}{1}$, is an isomorphism, and so $\ga_i =
\ass (S_i)$. We identify the rings $R_i$ and $S_i^{-1}R$ via the
isomorphism above. Since $S_i\in \maxDen_l(R)$ and $S_i^{-1}R$  is
a division ring, we see  that (by Lemma \ref{d1Dec12}) $R=S_i\coprod
\ga_i$, a disjoint union. Then
$$ S_i\cap\bigcap_{j\neq i} \ga_j = (\bigcap_{j\neq i}\ga_j)\backslash \{ 0\} \neq
\emptyset$$ since  $\bigcap_{i=1}^n\ga_i = \ker (\s ) =0$.  This
fact together with the fact that $R=S_i\coprod\ga_i$ for $i=1,
\ldots , n$,
 imply that the sets $S_1, \ldots , S_n$ are distinct.

 Finally, we claim that  $\maxDen_l(R) = \{ S_1, \ldots , S_n\}$. Suppose
that this is not true, that is there exists a set $S\in
\maxDen_l(R)$ distinct from the sets $S_1, \ldots , S_n$, we seek
a contradiction. For each element $i=1, \ldots , n$, by using the
disjoint union $R= S_i\coprod\ga_i$ and the fact that $S\neq S_i$,
we can find an element, say $t_i$ such that $t_i \in S\cap \ga_i$
(notice that $S\not\subseteq S_i$ since $S, S_i\in \maxDen_l(R)$).
Then $\s_i(t_i)=0$, and so $\s (t)=0$ where  $t:=t_1, \cdots
t_n\in S$. On the other hand, $\s (t)\neq 0$ since $t\in S$, a
contradiction. $\Box $

{\bf A left  localizable ring $R$ with $|\maxDen_l(R)|<\infty$ is a semiprime left Goldie
ring such that  $Q_{l,cl}(R)$ is a direct product of division rings, and vice versa}.
\begin{theorem}\label{3Dec12}
Let $R$ be a ring. The following statements are equivalent.
\begin{enumerate}
\item The ring $R$ is a left localizable ring with
$n:=|\maxDen_l(R)|<\infty$.\item$Q_{l,cl}(R)=R_1\times\cdots
\times R_n$ where $R_i$ are division rings.
 \item The ring $R$ is a semiprime left Goldie ring with $\udim
 (R)=|\Min (R)|=n$ where $\Min (R)$ is the set of minimal prime
 ideals of the ring $R$.
 \item $Q_l(R)= R_1\times\cdots \times R_n$ where $R_i$ are division rings.
\end{enumerate}
\end{theorem}

{\it Proof}. $(2\Leftrightarrow 3)$ It is known and easy to prove (use Goldie's Theorem).

$(1\Rightarrow 2)$ By Proposition \ref{A3Dec12},
$C'^{-1}R=R_1\times\cdots \times R_n$ where $R_i$ are division
rings. In particular, $C'^{-1}R$ is a semisimple (Artinian) ring.
Clearly, $C'\subseteq S_0(R)$ where $S_0(R)$ is the largest
regular left Ore set of the ring $R$. Then
$$ Q_l(R) = S_0(R)^{-1}R = Q_l(C'^{-1}R) = Q_l(R_1\times\cdots
\times R_n)=R_1\times\cdots \times R_n.$$ Since $Q_l(R)$ is a
semisimple ring, $Q_{l,cl}(R) = Q_l(R)$, by  (Corollary 2.10,
\cite{larglquot}).

$(2\Rightarrow 1)$ The map $\s :R\ra Q_{l,cl}(R)=R_1\times
\cdots\times R_n$, $r\mapsto \frac{r}{1}$, is a monomorphism. Any
division ring is a left localization maximal ring. By Proposition
\ref{A8Dec12}.(4), $\maxDen_l(R) = \{T_1, \ldots , T_n\}$ where
$T_i=\s^{-1} (R_1\times\cdots \times R_i^*\times\cdots \times
R_n)$ and $R_i \simeq T_i^{-1}R$. Since $R_1, \ldots , R_n$ are
division rings and $\s = \prod_{i=1}^n \s_i: R\ra \prod_{i=1}^n
T_i^{-1}R$, $r\mapsto (\s_1(r), \ldots , \s_n(r))$, is a
monomorphism where $\s_i :R\ra T_i^{_-1}R$, $r\mapsto
\frac{r}{1}$, the ring $R$ is a left localizable ring (Theorem \ref{29Nov12}).

$(3\Leftrightarrow 4)$ This equivalence is a particular case of a more general result (Corollary 2.10,
\cite{larglquot}): {\em Let $R$ be a ring. Then  $Q_l(R)$ is a semisimple ring iff $Q_{l, cl}(R)$ is a semisimple ring, and in this case} $Q_l(R)=Q_{l, cl}(R)$.  $\Box $


{\bf A bijection between $\maxDen_l(R)$ and $\maxDen_l(Q_l(R))$}.
\begin{proposition}\label{A8Dec12}
\cite{Crit-S-Simp-lQuot} Let $R$ be a ring, $S_l$ be the  largest regular left Ore set of the ring $R$, $Q_l:= S_l^{-1}R$ be the largest left quotient ring of the ring $R$, and $\CC$ be the set of regular elements of the ring $R$. Then
\begin{enumerate}
\item $S_l\subseteq S$ for all $S\in \maxDen_l(R)$. In particular,
$\CC\subseteq S$ for all $S\in  \maxDen_l(R)$ provided $\CC$ is a
left Ore set. \item Either $\maxDen_l(R) = \{ \CC \}$ or,
otherwise, $\CC\not\in\maxDen_l(R)$. \item The map $$
\maxDen_l(R)\ra \maxDen_l(Q_l), \;\; S\mapsto SQ_l^*=\{ c^{-1}s\,
| \, c\in S_l, s\in S\},
$$ is a bijection with the inverse $\CT \mapsto \s^{-1} (\CT )$
where $\s : R\ra Q_l$, $r\mapsto \frac{r}{1}$, and $SQ_l^*$ is the
sub-semigroup of $(Q_l, \cdot )$ generated by the set  $S$ and the
group $Q_l^*$ of units of the ring $Q_l$, and $S^{-1}R= (SQ_l^*)^{-1}Q_l$.
    \item  If $\CC$ is a left Ore set then the map $$ \maxDen_l(R)\ra \maxDen_l(Q), \;\; S\mapsto SQ^*=\{ c^{-1}s\,
| \, c\in \CC, s\in S\}, $$ is a bijection with the inverse $\CT
\mapsto \s^{-1} (\CT )$ where $\s : R\ra Q$, $r\mapsto
\frac{r}{1}$, and $SQ^*$ is the sub-semigroup of $(Q, \cdot )$
generated by the set  $S$ and the group $Q^*$ of units of the ring
$Q$, and $S^{-1}R= (SQ^*)^{-1}Q$.
\end{enumerate}
\end{proposition}





\begin{theorem}\label{C3Dec12}
Let $R$ be a ring, $T\in \Den_l(R, 0)$ and $\s : R\ra T^{-1}R$,
$r\mapsto \frac{r}{1}$. Then
\begin{enumerate}
\item $T\subseteq S$ for all $S\in \maxDen_l(R)$.  \item  The map
$$ \maxDen_l(R)\ra \maxDen_l(T^{-1}R), \;\; S\mapsto
\widetilde{S},$$ is a bijection with the inverse $\CT \mapsto
\s^{-1} (\CT )$ where $\widetilde{S}$ is the multiplicative monoid
generated in the ring $T^{-1}R$ by $\s (S)$ and $\s (T)^{-1}:=\{
t^{-1} \, | \, t\in T\}$ and $S^{-1}R\simeq \widetilde{S}^{-1}(T^{-1}R)$, and $S^{-1}R\simeq \widetilde{S}^{-1}(T^{-1}R)$.
\end{enumerate}
\end{theorem}

{\it Proof}. 1. Lemma \ref{1a27Nov12}.

2. Statement 2 follows from statement 1 and Proposition 3.4.(1),
\cite{larglquot}. $\Box $


\begin{theorem}\label{B3Dec12}
Let $R$ be a ring, $T\in \Den_l(R, 0)$ be such that $T^{-1}R=
R_1\times \cdots \times  R_n$ is a direct product of left
localization maximal rings; and $\tau : R\ra T^{-1}R$, $r\mapsto
\frac{r}{1}=(r_1, \ldots , r_n)$. Then $\maxDen_l(R) = \{ T_1,
\ldots , T_n\}$ where $T_i= \tau^{-1}(S_i)$, $S_i:=
R_1\times\cdots \times R_i^*\times\cdots \times R_n$ and
$T_i^{-1}R\simeq R_i$.
\end{theorem}

{\it Proof}. Proposition follows from Theorem \ref{C3Dec12} and
Theorem \ref{25Nov12}.(1). $\Box $

$\noindent $

The next theorem shows that $|\maxDen_l(R)|<\infty$ for all
semiprime left Goldie rings.

\begin{corollary}\label{a4Dec12}
Let $R$ be a semiprime left Goldie ring, $\CC$ be the set of
regular element of the ring $R$, $Q_{l,cl}=\prod_{i=1}^n Q_i$
where $Q_i$ are simple Artinian rings (by Goldie's Theorem), $\s
:R\ra Q_{l,cl}$, $r\mapsto \frac{r}{1}$. Then
\begin{enumerate}
\item $C\subseteq S$ for all $S\in \maxDen_l(R)$.\item
$\maxDen_l(R)= \{ S_1, \ldots , S_n\}$, $\maxDen_l(Q_{l,cl}(R))=
\{ S_1', \ldots , S_n'\}$ where $S_i=\s^{-1}(S_i')$ and
$S_i':=Q_1\times\cdots \times Q_i^*\times\cdots \times Q_n$ for
$i=1, \ldots , n$. The map
$$ \maxDen_l(R)\ra \maxDen_l(Q_{l,cl}(R)), \;\; S\mapsto
\widetilde{S},$$ is a bijection with inverse $\CT \mapsto \s^{-1}
(\CT )$ where $\widetilde{S}$ is the multiplicative monoid
generated in $Q_{l,cl}(R)$ by $\s (S)$ and $\s (\CC )^{-1}=\{
c^{-1}\, | \, c\in \CC\}$. \item $S_i^{-1}R\simeq Q_i$ for $i=1,
\ldots , n$.
\end{enumerate}
\end{corollary}

{\it Proof}.  This is a particular case of Theorem \ref{C3Dec12}
where $T=C$. $\Box $



\section{The core of a left Ore set}\label{TCLOS}

The aim of this section is to establish several properties of  the {\em core} $S_c$ of a left Ore set $S$ of a ring $R$, and when it exists (for example for a right Noetherian ring $R$) it is a powerful tool in studying localizations (Theorem \ref{A2Dec12}). A characterization of the core is given (Proposition \ref{B2Dec12}). The core shares remarkable properties (Lemma \ref{b2Dec12}, Theorem \ref{A2Dec12}).

$\noindent $

{\it Definition}, \cite{Crit-S-Simp-lQuot}. Let $R$ be a ring and $S$ be its left Ore set. The subset of $S$,
$$ S_c:=\{ s\in S\, | \, \ker (s\cdot ) = \ass (S)\}$$ is called the {\em core} of the left Ore set $S$ where $s\cdot : R\ra R$, $r\mapsto sr$.

$\noindent $

If $S_c\neq \emptyset$ then $S_cS_c\subseteq S_c$, i.e. the core $S_c$ is a multiplicative set. If $\ass (S)=0$ then $S_c=S$.

\begin{lemma}\label{b2Dec12}
If $S\in \Den_l(R)$ and $S_c\neq \emptyset$ then
\begin{enumerate}
\item $SS_c\subseteq S_c$.
\item For any $s\in S$ there exists an element $t\in S$ such that $ts \in S_c$.
\end{enumerate}
\end{lemma}

{\it Proof}. 1. Trivial.

2. Statement 2 follows directly from the left Ore condition: fix
an element $s_c\in S_c$, then $ts= rs_c\in S$ for some elements $t\in
S$ and $r\in R$. Since $\ass (S)\supseteq \ker ( ts\cdot
)=\ker(rs_c\cdot ) \supseteq \ker (s_c\cdot ) = \ass (S)$, i.e.
$\ker (ts\cdot ) = \ass (S)$, we have $ts \in S_c$.  $\Box $


\begin{theorem}\label{A2Dec12}
Suppose that $S\in \Den_l(R, \ga )$ and $S_c\neq \emptyset$. Then
\begin{enumerate}
\item $S_c\in \Den_l(R, \ga )$.
\item The map $\th : S_c^{-1}R\ra S^{-1}R$, $s^{-1}r\mapsto s^{-1}r$, is a ring $R$-isomorphism. So, $S_c^{-1}R\simeq S^{-1}R$.
\end{enumerate}
\end{theorem}

{\it Proof}. 1. By Lemma \ref{b2Dec12}.(1), $S_cS_c\subseteq S_c$,
that is the set $S_c$ is a multiplicative set. By Lemma
\ref{b2Dec12}.(2), $S_c\in \Ore_l(R)$: for any elements $s_c\in
S_c$ and $r\in R$, there are elements $s\in S$ and $r'\in R$ such
that $sr = r's_c$ (since $S\in \Ore_l(R)$). By Lemma
\ref{b2Dec12}.(2), $s_c':= ts\in S_c$ for some $t\in S$. Then $s_c' r=
tr' s_c$.

If $rs_c=0$ for some elements $r\in R$ and $s_c\in S_c$ then $sr=0$ for some element $s\in S$ (since $S\in \Den_l(R)$). By Lemma \ref{b2Dec12}.(2), $s_c':=ts\in S_c$ for some element $t\in S$, hence $s_c'r=0$. Therefore, $S_c\in \Den_l(R, \ga )$.

2. By statement 1 and the universal property of left Ore
localization,  the map $\th$ is a well-defined monomorphism. By
Lemma \ref{b2Dec12}.(2), $\th$ is also a surjection: let
$s^{-1}r\in S^{-1}R$, and $s_c:= ts\in S_c$ for some element $t\in
S$ (Lemma \ref{b2Dec12}.(2)). Then
$$ s^{-1}r= s^{-1}t^{-1}tr= (ts)^{-1}tr= s_c^{-1} tr. \;\; \Box $$


Let $S$ be a left Ore set of a ring $R$. The set of  right annihilators $\ker (S\cdot ) :=\{ \ker (s\cdot ) \, | \, s\in S\}$ is a poset with respect to $\subseteq $. Let
$$ \max (S) :=\{ s\in S\, | \, \ker (s\cdot ) \;\; {\rm is \; a\; maximal \;element \; of}\;\; \ker (S\cdot )\}.$$
The set $\max (S)$ is a non-empty set for a ring with ACC on right annihilators (for example, when $R$ is a right Noetherian ring).

\begin{lemma}\label{a2Dec12}
Let $S\in \Ore_l(R)$. Then $\ass (S) = \bigcup_{s\in S}\ker (s\cdot ) = \sum_{s\in S}\ker (s\cdot )$.
\end{lemma}

{\it Proof}. Let $U$ and $\S$ stand for the union and the sum respectively. Then $\ass (S) \supseteq \S \supseteq U= \ass (S)$, and so $\ass (S) = \S = U$. $\Box $


\begin{proposition}\label{B2Dec12}
Let $S$ be a left Ore set of a ring $R$. Then $S_c= \max (S)$.
\end{proposition}

{\it Proof}. If $S_c\neq \emptyset$ then $S_c = \max (S)$, by Lemma \ref{a2Dec12}.

Suppose that $S_c = \emptyset$. To finish the proof we have to show that $\max (S) = \emptyset$. Suppose that this is not
 the case, we seek a contradiction. Let $\ga = \ass (S)$. Fix $s\in \max (S)$. Then $\gb := \ker (s\cdot ) \subseteq \ga$.
 We claim that $\gb = \ga $, i.e. $s\in S_c$, a contradiction. To prove the claim we have to show that for any element
 $t\in S$, $\ker (t\cdot ) \subseteq \gb$ (see Lemma \ref{a2Dec12}). Notice that $\ker (s's\cdot ) = \gb$ for all $s'\in S$,
 by the maximality of $\ker (s\cdot )$. Using the left Ore condition for $S$, we have the equality $s_1s= rt$ for some elements $s_1\in S$ and $r\in R$. Hence,
$$ \gb = \ker (s_1s\cdot ) = \ker (rt\cdot ) \supseteq \ker (t\cdot ), $$ as required.  $\Box $


The next theorem gives an explicit description of the cores of maximal left denominator sets of a left localizable ring provided $|\maxDen_l(R)|<\infty$.

\begin{theorem}\label{C2Dec12}
Let $R$ be a left localizable ring such that $\maxDen_l(R) = \{ S_1, \ldots , S_n\}$.
\begin{enumerate}
\item If $n=1$ then $S_{1,c} = S_1= R\backslash \{ 0\}$.
\item If $n\geq 2$ then $S_{i, c} = S_i\cap \bigcap_{j\neq i} \ga_i$  where $\ga_j=\ass (S_j)$.
\end{enumerate}
\end{theorem}

{\it Proof}. We keep the notation of Theorem \ref{29Nov12} and its proof.

1. If $n=1$ then $R$ is a domain $S_1= R\backslash \{ 0\}$
(Theorem \ref{29Nov12}), and so $S_{1, c}=  S_1=R\backslash \{
0\}$.

2. Suppose that $n\geq 2$. For each $i=1, \ldots , n$, $C_i':=
S_i\cap \bigcap_{j\neq i} \ga_j\neq \emptyset$, by Corollary
\ref{e1Dec12}. By Proposition \ref{A3Dec12}.(2), $C_i'\subseteq
S_{i,c}$.  By Theorem \ref{29Nov12}, the map
$$ \s : =\prod_{i=1}^n \s_i : R\ra \prod_{i=1}^n R_i , \;\; r\mapsto (r_1, \ldots , r_n), $$ is a ring monomorphism.
 Since $R_i=S_i^{-1}R$ is a division ring (Theorem \ref{29Nov12}), $R=S_i\coprod \ga_i$, by Lemma \ref{d1Dec12}.
 Clearly, $\ga_i =\{ r= (r_1, \ldots , r_n) \, | \, r_i=0\}$. Each element $s'$ of the set $\CC_i'$ has the form
 $(0, \ldots , 0 , s_i', 0, \ldots , 0)$ with $s_i'\in R_i^*$. Then clearly, $s'\ga_i=0$ and so $C_i'\subseteq S_{i, c}$.
  To show that the equality $S_i'= S_{i, c}$ holds it suffices to show that every element $s\in S_i\backslash C_i'$ does
  not not belong to $S_{i, c}$. Fix $s$ such that $s\in S_i\backslash C_i'$. Then there is an index, say  $j$,  such that $j\neq i$ and
  such that $s_j\neq 0$ in $s= (s_1, \ldots , s_n)$. Then $s\cdot \CC_j' \neq 0$ (Corollary \ref{e1Dec12}) but $\CC_j'\subseteq \ga_i$. Therefore, the element $s$
  does not belong to $S_{i, c}$.  $\Box $


$${\bf Acknowledgements}$$

 The work is partly supported by  the Royal Society  and EPSRC.


\small{

Department of Pure Mathematics

University of Sheffield

Hicks Building

Sheffield S3 7RH

UK

email: v.bavula@sheffield.ac.uk}

\end{document}